\documentclass{article}
\PassOptionsToPackage{dvipsnames}{xcolor}
\usepackage{graphicx} 

\usepackage{amsmath, amsthm, amssymb, amsfonts}
\usepackage{thmtools}
\usepackage{graphicx}
\usepackage{setspace}
\usepackage{appendix}
\usepackage{geometry}
\usepackage[mathscr]{euscript}
\usepackage{float}
\usepackage{tikz-cd}
\usepackage{hyperref}
\usepackage[utf8]{inputenc}
\usepackage[english]{babel}
\usepackage{framed}
\usepackage{relsize}
\usepackage{tcolorbox}
\usepackage{xcolor}
\usepackage{indentfirst}
\usepackage{dynkin-diagrams}
\usepackage{tikz}
\usepackage{tikzsymbols}
\usepackage[backend=biber,style=numeric]{biblatex}
\newcommand\C{\ensuremath{\mathbb{C}}}

\newcommand\msc{\mathscr{C}}

\newcommand\Spec{\operatorname{Spec}}

\newcommand\colim{\lim\limits_{\ar}}

\newcommand\mco{\mathcal{O}}

\newcommand\ar{\longrightarrow}

\newcommand\eps{\epsilon}

\newcommand\mf{\mathfrak}
\newcommand\Z{\ensuremath{\mathbb{Z}}}

\newcommand\g{\mathfrak{g}}

\newcommand\Hom{\text{Hom}}
\newcommand\h{\mathfrak{h}}
\newcommand\st{\text{ }|\text{ }}
\newcommand{\DynkinBn}[1]

\title{\textbf{An illustration of formal moduli problems with differential graded Lie algebras}}
\author{ ethan eugene wynner}

\begin{document}

\maketitle 

Abstract: This article provides an exposition to the topic of formal moduli problems, emphasizing its connections with differential graded Lie algebras, and mainly following from Jacob Lurie's \textit{DAG X: Formal Moduli Problems}. As such, this paper should not be viewed as a presentation of original work, but rather a concise introduction to the subject in the form of a set of organized notes. I hoped to make this paper feel welcoming and insightful for the non-expert enjoyer of derived algebraic geometry, like myself. Enjoy! \Laughey  

\newpage
\begin{center}
    \textit{"In} 1984, \textit{I was hospitalized for approaching perfection."}\newline\newline -David Berman, \textit{Random Rules}
\end{center}
\newpage 

\tableofcontents

\newpage
\section{Deformation theory and formal moduli problems}
First we note a down to earth characterization of the study of deformations, and leap into the more ethereal realm of deformation theory connected to the study of formal moduli problems.
\subsection{First order deformations and the ring of dual numbers}
Let $k$ be a field. The most straightfoward-to-the-point object of study in deformation theory is the ring of dual numbers $k[t]/(t^2).$ The idea is to suppose that we are given some algebro-geometric structure (i.e., a scheme). Our deformation-theoretic impulse then is to try to classify extensions of this structure over $k[t]/(t^2).$ These extensions are called \textit{first order deformations.} A technical consideration is that we must assume that the extensions are \textit{flat} over $k[t]/(t^2).$ We will define this now. 

\subsubsection{Definition [5]} 
Let $A$ be a ring. We call an $A$-module $M$ \textit{flat} if the functor 
$$N \longmapsto N \otimes_A M$$
is exact on the category of $A$-modules. We say that a morphism $f: X \ar Y$ of schemes is flat if for every point $x \in X,$ the local ring $\mco_{X,x}$ is flat over the local ring $\mco_{Y, f(x)}.$ Finally, a sheaf of $\mco_X$-modules $\mathscr{F}$ is flat over $Y$ if for every $x \in X,$ the stalk $\mathscr{F}_x$ is flat over $\mco_{Y,f(x)}.$

\subsubsection{Example-definition of a deformation (problem)}
Let $X$ be a scheme over $k,$ and let $Y$ be a closed subscheme of $X.$ A \textit{deformation of $Y$} over $k[t]/(t^2)$ in $X$ is a closed subscheme $Y \subset X \times k[t]/(t^2).$ The "problem" here is to classify all deformations of $Y$ over $k[t]/(t^2).$

\subsection{The Rocketship}
We are now interested in approaching deformation theory from the standpoint of higher category theory. In so doing, we still wish to study classes of algebro-geometric objects, but, since we are so infatuated with Grothendieck, we identify these objects with functors
$$X: \Gamma \ar \mathscr{S},$$
where $\Gamma$ is some presentable $\infty$-category of these aforementioned objects, and $\mathscr{S}$ is the $\infty$-category of spaces (see appendix B). We ask little of the category $\Gamma;$ only that it has suitable objects and a terminal object $*.$ This is so as to enable us to view a \textit{point} of a functor $X$ to be a point in the space $X(*).$ We also denote by $\Gamma_*$ the $\infty$-category of pointed objects in $\Gamma,$ that is, pairs $(Y, \eta)$ where $Y$ is an object of $\Gamma$ and $\eta: * \ar Y$ is the unique map. We also take into high consideration the forgetful functors $\Omega^{\infty-n}_*: \text{Stab}(\Gamma) \ar \Gamma_*$ defined on all $n \in \Z.$ Moreover, we have $\Omega^{\infty - n}: \text{Stab}(\Gamma) \ar \Gamma$ given by composition with the forgetful functor $\Gamma_* \ar \Gamma.$ This brings us to a trampoline which we will use to attain the idea of a formal moduli problem: 

\subsubsection{Definition}
A \textit{deformation context} is a pair $(\Gamma,\{E_\alpha\}_{\alpha \in T})$ where $\Gamma$ follows our previous rules, and $\{E_\alpha\}_{\alpha \in T}$ is a set of objects in $\text{Stab}(\Gamma).$
\newline\newline 
Let $(\Gamma,\{E_\alpha\}_{\alpha \in T})$ be a deformation context. A morphism $A' \ar A$ in $\Gamma$ is called \textit{elementary} if there exists $\alpha \in T, n> 0,$ and a pullback diagram
$$\begin{tikzcd}
	{A'} & {*} \\
	A & {\Omega^{\infty-n}E_\alpha}
	\arrow[from=1-1, to=1-2]
	\arrow["\phi", from=1-1, to=2-1]
	\arrow["{\phi_0}", from=1-2, to=2-2]
	\arrow[from=2-1, to=2-2]
\end{tikzcd}$$

\subsubsection{Observation} Let $A \in \Gamma.$ Every elementary map $A' \ar A$ in $\Gamma$ is given by the fiber of a map $A \ar \Omega^{\infty-n}E_\alpha$ for some $n> 0, \alpha \in T.$
\subsubsection{Definition}
We say that a morphism $\phi: A' \ar A$ is \textit{small} if it can be written as a finite composition of elementary morphisms $A'\simeq A_0 \ar A_1 \ar \cdots \ar A_n \simeq A.$ We say that $A \in \Gamma$ is small (as an object of $\Gamma$) if the map $A \ar *$ is small. We denote by $\Gamma^{sm}$ the full subcategory of $\Gamma$ spanned by small objects. 
\newline\newline 
We now give a definition of the titular object of this paper:

\subsubsection{A general interpretation of a formal moduli problem}
Let $(\Gamma,\{E_\alpha\}_{\alpha \in T})$ be a deformation context. A \textit{formal moduli problem} is a functor
$$X: \Gamma^{sm} \ar \mathscr{S}$$
satisfying 
\begin{enumerate}
    \item The space $X(*)$ is contractible
    \item If $$\begin{tikzcd}
	{A'} & {B'} \\
	A & B
	\arrow[from=1-1, to=1-2]
	\arrow[from=1-1, to=2-1]
	\arrow["\phi", from=1-2, to=2-2]
	\arrow[from=2-1, to=2-2]
\end{tikzcd}$$ is a pullback diagram and $\phi: B' \ar B$ is small, then the image 
$$\begin{tikzcd}
	{X(A')} & {X(B')} \\
	{X(A)} & {X(B)}
	\arrow[from=1-1, to=1-2]
	\arrow[from=1-1, to=2-1]
	\arrow["{X(\phi)}", from=1-2, to=2-2]
	\arrow[from=2-1, to=2-2]
\end{tikzcd}$$ is a pullback diagram in $\mathscr{S}.$
\end{enumerate}
Let $\text{Moduli}^\Gamma$ denote the full subcategory of $\text{Fun}(\Gamma^{sm}, \mathscr{S})$ spanned by formal moduli problems. We call $\text{Moduli}^\Gamma$ the \textit{$\infty$-category of formal moduli problems.}

\subsubsection{Example}
Let $(\Gamma,\{E_\alpha\}_{\alpha \in T})$ be a deformation context, $A \in \Gamma$ an object. Let $\Spec A: \Gamma \ar \mathscr{S}$ be the functor corepresented by $A,$ given on $\Gamma^{sm}$ by $\Spec A(B) = \text{Map}_\Gamma(A,B).$ Then $\Spec A$ is a formal moduli problem, and, moreover, the assignment $A \longmapsto \Spec A$ determines a functor $\Spec: \Gamma^{op} \ar \text{Moduli}^\Gamma.$

\subsubsection{Admissory note:}
For much of what we wish to do in this paper, it will be necessary to appeal to notions of spectral algebraic geometry, such as $\mathbb{E}_\infty$-rings/algebras; an example of this being the $\infty$-category $\text{CAlg}^{aug}_k$ of augmented $\mathbb{E}_\infty$-algebras over an $\mathbb{E}_\infty$-ring $k,$ as well as the $\infty$-category $\text{Mod}_k$ of $k$-module spectra, both of which are crucial to some of the important stuff in here.   A full and proper (in the non-mathematical sense) exposition of this material is literally unfeasible given the focus and length of this paper, although I have tried to make some reasonable accommodations when necessary. Fortunately, there is an abundance of resources, as arrayed in the references section.

\subsection{The tangent complex}
We will now delve into a very important framework upon which we can rest the laurels of formal moduli problems. In order to do so, we will generalize the construction of the \textit{Zariski tangent space,} which we will first review.
\subsubsection{The Zariski tangent space}
Let $X$ be an algebraic variety over $\C,$ and let $x: \Spec \C \ar X$ be a point of $X.$ A tangent vector to $X$ at $x$ is a dashed arrow making the diagram 
$$\begin{tikzcd}
	{\Spec \C} & X \\
	{\Spec \C[t]/(t^2)} & {\Spec \C}
	\arrow[from=1-1, to=1-2]
	\arrow[from=1-1, to=2-1]
	\arrow[from=1-2, to=2-2]
	\arrow[dashed, from=2-1, to=1-2]
	\arrow[from=2-1, to=2-2]
\end{tikzcd}$$
commute. The collection of tangent vectors at $x$, denoted $T_xX,$ is called the \textit{Zariski tangent space.} Moreover, if $\mco_{X,x}$ is the local ring of $X$ at $x,$ with $\mf{m}$ its maximal ideal, there is a canonical isomorphism $T_xX \simeq (\mf{m}/\mf{m}^2)^\vee,$ where $(\mf{m}/\mf{m}^2)^\vee$ is called the \textit{cotangent space} to $X.$ Our first step in generalizing this construction into the setting of formal moduli problems is to associate the variety $X$ with its 'functor of points' $X(A) = \Hom_{\text{Sch}(\C)}(\Spec A, X),$ where $\text{Sch}(\C)$ is the category of schemes over $\C.$ Then $T_xX$ can be characterized as the fiber of the map 
$$X(\C[t]/(t^2)) \ar X(\C)$$
over the point $x \in X(\C).$ Note also that in a mercurial way the (commutative) ring of dual numbers $\C[t]/(t^2)$ is given by $\Omega^\infty E,$ where $E$ is the spectrum object in $\text{CAlg}^{aug}_\C$ corresponding to $\C$ in the sense that the suspension functor defined by the pushout $\Sigma\C = \C\coprod_\C\C \simeq \C[t]/(t^2)$, where we take $E = \Sigma^\infty(\C) \in \text{Stab}(\text{CAlg}^{aug}_\C)$ so that $\Omega^{\infty-0}(E) = \Sigma\C \simeq \C[t]/(t^2).$
This determines a deformation context $(\text{CAlg}^{aug}_\C, \{E\}),$ and leads us to suspect a generalization.

\subsubsection{Generalization}
Let $(\Gamma,\{E_\alpha\})$ be a deformation context. Let $Y: \Gamma^{sm} \ar \mathscr{S}$ be a formal moduli problem. For every $\alpha,$ the (generalized) \textit{tangent space} of $Y$ at $\alpha$ is the space $Y(\Omega^\infty E_\alpha).$
\subsubsection{Definition}
Let $\msc$ be an $\infty$-category with finite colimits, and $\mathscr{D}$ an $\infty$-category with finite limits. We say that a functor $F: \msc \ar \mathscr{D}$ is \textit{excisive} if for every pushout diagram $P$ in $\msc,$ the image $F(P)$ is a pushout diagram in $\mathscr{D}.$ We say that $F$ is \textit{strongly excisive} if it is excisive and maps initial objects to terminal objects.
\newline\newline Let $\mathscr{S}^{fin}_*$ be the $\infty$-category of finite pointed sets. If $\mathscr{D}$ is an $\infty$-category with finite limits, then $\text{Stab}(\mathscr{D})$ is the full subcategory of $\text{Fun}(\mathscr{S}^{fin}_*, \mathscr{D})$ spanned by pointed excisive functors. In particular, we can identify the $\infty$-category \textbf{Spc} = $\text{Stab}(\mathscr{S})$ of spectra in this way: as the full subcategory of $\text{Fun}(\mathscr{S}^{fin}_*, \mathscr{S})$ spanned by strongly excisive functors.

\subsubsection{Construction}
Let $(\Gamma, \{E_\alpha\})$ be a deformation context. For all $\alpha,$ we can identify $E_\alpha \in \text{Stab}(\Gamma)$ with the functor $E_\alpha: \mathscr{S}^{fin}_*, \Gamma.$ Then 
\begin{enumerate}
    \item For any map $f: K' \ar K'$ of finite pointed spaces inducing a surjection on homotopy groups $\pi_0K \ar \pi_0K',$ the induced map $E_\alpha(K) \ar E_\alpha(K')$ is small in $\Gamma.$
    \item For all $K \in \mathscr{S}^{fin}_*, E_\alpha(K) \in \Gamma$ is small.
\end{enumerate}

\subsubsection{Definition}
Let $\Gamma^{sm} \ar \mathscr{S}$ be a formal moduli problem. For every $\alpha,$ we view the composition 
$$Y(E_\alpha)= \mathscr{S}^{fin}_* \xrightarrow{E_\alpha}\Gamma^{sm} \xrightarrow{Y} \mathscr{S} $$
as an object of $\textbf{Spc},$ and we call $Y(E_\alpha)$ the \textit{tangent complex} to $Y$ at $\alpha.$

\subsubsection{Remark}
We can further identify the tangent space $Y(\Omega^\infty E_\alpha)$ with the 0th rung of the tangent complex $Y(E_\alpha).$ Moreover, there are pleasing canonical homotopy equivalences
$$Y(\Omega^{\infty-n}E_\alpha) \simeq \Omega^{\infty - n}Y(E_\alpha).$$

\subsection{Deformation Theories}
We will introduce the notion here of (weak) deformation theories, as well as state a promising theorem. Our setup is as follows. Let $(\Gamma, \{E_\alpha\})$ be a deformation context. Suppose that $\Xi$ is an $\infty$-category. We want to figure out when there is an equivalence $\text{Moduli}^\Gamma \simeq \Xi.$ To each $A \in \Gamma,$ we associate the formal moduli problem $\Spec A \in \text{Moduli}^\Gamma,$ given by $\Spec A(R) = \text{Map}_\Gamma(A,R).$ Moreover, if $\text{Moduli}^\Gamma \simeq \Xi,$ we obtain a functor 
$$\mathfrak{D}: \Gamma^{op} \ar \Xi.$$

\subsubsection{Definition}
Let $(\Gamma, \{E_\alpha\})$ be a deformation context. A \textit{weak deformation theory} is a functor $\mf{D}: \Gamma^{op} \ar \Xi $ satisfying the following conditions:
\begin{enumerate}
    \item $\Xi$ is presentable
    \item $\mf{D}$ admits a left adjoint $\mf{D}': \Xi \ar \Gamma^{op}$
    \item $\Xi$ has a full subcategory $\Xi_0$ such that 
    \begin{enumerate}
        \item For each $K \in \Xi_0,$ the unit map $K \longmapsto \mf{D}(\mf{D}'(K))$ is an equivalence.
        \item The initial object $\varnothing \in \Xi_0.$ Hence $\varnothing \simeq \mf{D}(\mf{D}'(\varnothing)) \simeq \mf{D}(*).$
        \item For all $\alpha, n \ge 1,$ there exists an object $K_{\alpha,n} \in \Xi$ and an equivalence $\Omega^{\infty - n}E_\alpha\simeq \mf{D}'K_{\alpha,n},$ determining a map 
        $$v_{\alpha,n}: K_{\alpha, n} \simeq \mf{D}(\mf{D}'(K_{\alpha,n})) \simeq D(\Omega^{\infty - n}E_\alpha) \ar D(*) \simeq \varnothing.$$
        \item For every pushout 
        $$\begin{tikzcd}
	{K_{\alpha,n}} & K \\
	\varnothing & {K'}
	\arrow[from=1-1, to=1-2]
	\arrow["{v_{\alpha,n}}", from=1-1, to=2-1]
	\arrow[from=1-2, to=2-2]
	\arrow[from=2-1, to=2-2]
\end{tikzcd}$$
if $K \in \Xi_0,$ then $K' \in \Xi_0.$
    \end{enumerate}
\end{enumerate}

\subsubsection{Example} Let $k$ be a field of characteristic zero and $(\text{CAlg}^{aug}_k, \{E\})$ be the deformation context from earlier. We will later construct a weak deformation theory 
$$\mf{D}: (\text{CAlg}^{aug}_k)^{op} \ar \text{Lie}_k$$
(where $\text{Lie}_k$ is the $\infty$-category of differential graded Lie algebras), where the adjoint is the cohomological Chevalley Eilenberg functor. In fact, this assignment $\g_* \ar C^*(\g_*)$ is actually a deformation theory, as it satisfies an extra condition which we will get to in a moment.

\subsubsection{Moreover-proposition}
Let $\mf{D}:\Gamma^{op} \ar \Xi$ be a weak deformation theory. Condition 3 above implies the following
\begin{enumerate}
    \item $\mf{D}$ carries terminal objects in $\Gamma$ to initial objects in $\Xi.$
    \item Let $A = \mf{D}'(K) \in \Gamma, K \in \Xi.$ The unit map $\mf{D}'(\mf{D}(A))$ is an eqiuivalence in $\Gamma.$
    \item If $A \in \Gamma$ is small, then $\mf{D}(A) \in \Xi_0,$ and $A \ar \mf{D}'(\mf{D}(A))$ is an equivalence in $\Gamma.$
    \item If $\sigma = $
    $$\begin{tikzcd}
	{A'} & {B'} \\
	A & B
	\arrow[from=1-1, to=1-2]
	\arrow[from=1-1, to=2-1]
	\arrow["\phi", from=1-2, to=2-2]
	\arrow[from=2-1, to=2-2]
\end{tikzcd}$$
is a pullback diagram where $A,B, \phi$ are small, then $\mf{D}(\sigma)$ is a pushout in $\Xi.$
\end{enumerate}
\subsubsection{A pair of corollaries}
\begin{enumerate}
    \item Let $y: \Xi \ar \text{Fun}(\Xi,\mathscr{S})$ be the Yoneda embedding. For every $K \in \Xi,$ the composition
    $$\Gamma^{sm} \subset \Gamma \xrightarrow{\mf{D}} \Xi^{op} \xrightarrow{y(K)} \mathscr{S}$$
    is a formal moduli problem, which determines a functor $\Psi: \Xi \ar \text{Moduli}^\Gamma \subset \text{Fun}(\Gamma^{op},\mathscr{S}). $
    \item Let $\mf{D}: \Gamma^{op} \ar \Xi$ be a weak deformation theory. For every $\alpha, K \in \Xi,$ the composition 
    $$\mathscr{S}^{fin}_* \xrightarrow{E_\alpha}\Gamma \xrightarrow{\mf{D}}\Xi^{op} \xrightarrow{y(K)}\mathscr{S}$$
    is strongly excisive, and can be indeitifed with a spectrum object $e_\alpha(K) \in \textbf{Spc}.$ This determines a functor 
    $$e_\alpha: \Xi \ar \textbf{Spc}.$$
\end{enumerate}

\subsubsection{Definition}
Finally we get here: A \textit{deformation theory} is a weak deformation theory satisfying one extra condition:
\newline\newline For every $\alpha,$ if $e_\alpha: \Xi \ar \textbf{Spc}$ is the functor above, then $e_\alpha$ preserves small sifted colimits, and a morphism $f: A \ar B \in \Xi$ is an equivalence iff $e_\alpha(f)$ is an equivalence; i.e. $e_\alpha(f): e_\alpha(A) \ar e_\alpha(B)$ is a weak equivalence in textbf{Spc}. This construction allows us to view \textbf{Spc} = $\text{Stab}(\mathscr{S}) \subset \text{Fun}_{\text{exc}}(\mathscr{S}^{fin}_*, \textbf{Spc})$ of excisive functors. 

\subsubsection{Theorem (Lurie)}
\textit{Let $\mf{D}: \Gamma^{op} \ar \Xi$ be a deformation theory. Then the functor} 
$$\Psi: \Xi \ar \text{Moduli}^\Gamma$$
\textit{is an equivalence of $\infty$-categories.}
\newline\newline The proof of this theorem can be found in Lurie [7] \S1.5

\subsubsection{Remark}
The composition
$$\Gamma^{op}\xrightarrow{\mf{D}}\Xi \xrightarrow{\Psi} \text{Moduli}^\Gamma$$
carries an object $A \in \Gamma$ to the formal moduli problem given by
$$B \longmapsto \text{Maps}_\Xi(\mf{D}(B), \mf{D}(A)) \simeq \text{Map}_\Gamma(A, \mf{D}'(\mf{D}(B)).$$
Hence the unit map $B \ar \mf{D}'(\mf{D}(B))$ determines a natural transformation $\beta: \Spec \ar \Psi \circ \mf{D}.$ It follows from proposition 1.4.3 that $\beta$ is an equivalence. Combined with theorem 1.4.6, we observe that $\mf{D}$ is equivalent to the weak deformation theory $\Spec: \Gamma^{op} \ar \text{Moduli}^\Gamma.$

\newpage

\section{Formal moduli problems for commutative algebras }
 \subsection{Introduction} Our main objective here is to connect the theory of formal moduli problems with that of differential graded Lie algebras. We offer the following proposition of great importance (PGI), which we hope to properly unwind over the course of this paper:

 \subsubsection{Proposition of great importance} 
 (PGI): \textit{If $X$ is a moduli space over a field $k$ of characteristic zero, then a formal neighborhood of any point $x \in X$ is controlled by a differential graded Lie algebra.}
\newline\newline
First of all, what do we mean by moduli space? For an example, let $k = \C.$ We can take $X$ to be a scheme, with a functor $R \longmapsto X(R) = \Hom(\Spec R, X),$ where $R$ is a commutative ring. We define a \textit{classical moduli problem} to be a functor 
$$X : \text{Ring}_\C \ar \textbf{Set}$$
where $\text{Ring}_\C$ is the category of commutative $\C$-algebras. For our purposes here, it is sometimes the case that a functor taking values in \textbf{Set} will not be adequate. As such, we define the following variant, which captures far more possibilities: 

\subsubsection{Definition}
Let $\msc$ be an $\infty$-category. A \textit{$\msc$-valued classical moduli problem} is a functor 
$$N(\text{Ring}_\C) \ar \msc,$$
where $N(\text{Ring}_\C)$ denotes the \textit{nerve} of the category $\text{Ring}_\C.$ We now see that our original definition is a special case of this new one, taking $\msc = N(\textbf{Set}).$ Now one might ask: what do we mean by a formal neighborhood? A pertinent and ubiquitous example follows: Let $k = \C,$ and let $X = \Spec A$ be an affine variety over $\C.$ A closed point $x \in X$ is determined by a $\C$-algebra homomorphism $\phi: A \ar \C,$ which itself determined by the choice of maximal ideal $\mf{m} = \ker(\phi) \subset A.$ The \textit{formal completion} of $X$ at the point $x$ is the functor $X^\wedge: \text{Ring} \ar \textbf{Set}$ given by taking $X^\wedge(R)$ to be the collection of commuting ring homomorphisms $A \ar R$ which carry elements of $\mf{m}$ to nilpotents in $R.$ That is, 
$$X^\wedge(R) = \{f \in X(R) \st f(\Spec R) \subset \{x\} \subset \Spec A\}.$$
\subsubsection{Definition} Let $R \in \text{Ring}_\C.$ We say that $R$ is \textit{local artinian} if it is finite dimensional as a $\C$-vector space, and is a local ring, i.e. has a unique maximal ideal $\mf{m}_R.$ We denote by $\text{Ring}^{art}_\C$ the category of local Artinian $\C$-algebras. Importantly, we observe that if $X$ is an affine variety over $\C,$ its formal completion $X^\wedge$ at $x \in X$ can be recovered by its values on local artinian rings. As such, we can further refine our definiton:

\subsubsection{Refined definition}
Let $\msc$ be an $\infty$-category. A \textit{$\msc$-valued classical formal moduli problem is a functor}
$$N(\text{Ring}^{art}_\C) \ar \msc.$$
If $X$ is \textbf{Set}-valued, and we have a point $\eta \in X(\C)$, we can define a \textbf{Set}-valued classical formal moduli problem $X^\wedge$ by 
$$X^\wedge(R) = X(R) \times_{X(R/\mf{m}_R)}\{\eta\},$$
which we call the \textit{completion of $X$ at $\eta.$}
Similarly, if $X$ is Gpd-valued (where Gpd denotes the category of groupoids), we can use the same formula using a homotopy fiber product.

\subsubsection{A palatable example}
Let $R$ be a commutative $\C$-algebra. Let $X(R)$ be the groupoid whose objects are smooth proper $R$-schemes and whose isomorphisms are those of such $R$-schemes. Take $\eta \in X(\C),$ corresponding to a smooth proper algebraic variety $Z.$ The functor $X^\wedge$ assigns to each $R \in \text{Ring}^{art}_\C$ the groupoid $X^\wedge(R)$ of \textit{deformations} over $Z$ (over $R$). That is, smooth proper morphisms $f: \overline{Z}\ar \Spec R$ fitting into the pullback diagram 
$$\begin{tikzcd}
	Z & {\overline{Z}} \\
	{\Spec \C} & {\Spec R}
	\arrow[from=1-1, to=1-2]
	\arrow[from=1-1, to=2-1]
	\arrow[from=1-2, to=2-2]
	\arrow[from=2-1, to=2-2]
\end{tikzcd}.$$
 
The functor $X^\wedge$ has some important properties:
\begin{enumerate}
    \item The image under $X^\wedge$ of the ring of dual numbers, $X^\wedge(\C[t]/(t^2),$ is the groupoid of first order deformations of the variety $Z.$ Moreover, every first order deformation $Z$ has an automorphism group which is naturally isomorphic to $H^0(Z; T_Z),$ where $T_Z$ is the tangent bundle of $Z.$
    \item The collection of isomorphism classes of first order deformations of $Z$ are naturally identified with the first cohomology $H^1(Z; T_Z).$
    \item Every first order deformation $\eta_1$ of $Z$ can be assigned a class $\theta \in H^2(Z; T_Z)$ which vanishes if and only if $\eta_1$ extends to a second order deformation $\eta_2  \in X^\wedge(\C[t]/(t^3)).$
\end{enumerate}
The first two of these properties are nice and friendly, and can be exposited without too much extra machinery. However, property 3 is rather unfriendly; in order to properly explain it, one must turn to the workings of spectral algebraic geometry, specifically the theory of commutative ring spectra, or $\mathbb{E}_\infty$-rings/ring spaces. Unfortunately for Refined definition 2.1.4, our construction is not complete for arbitrary classical formal moduli problems; we cannot assume $X^\wedge$ is defined on non-discrete $\mathbb{E}_\infty$-rings). This leads us to an even-more-refined definition, that of a \textit{formal moduli problem,} after which this paper is jointly named. 

\subsubsection{Definition of a formal moduli problem for commutative algebras}
Let $\text{CAlg}^{sm}_\C$ denote the category of small $\mathbb{E}_\infty$-algebras over $\C$. Let $\mathscr{S}$ denote the $\infty$-category of spaces. A \textit{formal moduli problem} over $\C$ is a functor 
$$X: \text{CAlg}^{sm}_\C \ar \mathscr{S}$$
satisfying the following two properties:
\begin{enumerate}
\item The space $X(\C)$ is contractible.
    \item For every pullback diagram 
    $$\begin{tikzcd}
	R & {R_0} \\
	{R_1} & {R_{01}}
	\arrow[from=1-1, to=1-2]
	\arrow[from=1-1, to=2-1]
	\arrow[from=1-2, to=2-2]
	\arrow[from=2-1, to=2-2]
\end{tikzcd}$$
in $\text{CAlg}^{sm}_\C$ for which the underlying maps $\pi_0R_0 \ar \pi_0R_{01}$ and $\pi_0R_1 \ar \pi_0R_{01}$ are surjective, the diagram 
$$\begin{tikzcd}
	{X(R)} & {X(R_0)} \\
	{X(R_1)} & {X(R_{01})}
	\arrow[from=1-1, to=1-2]
	\arrow[from=1-1, to=2-1]
	\arrow[from=1-2, to=2-2]
	\arrow[from=2-1, to=2-2]
\end{tikzcd}$$
admits a unique factorization
$$\begin{tikzcd}
	S \\
	& {X(R)} & {X(R_0)} \\
	& {X(R_1)} & {X(R_{01})}
	\arrow[dashed, from=1-1, to=2-2]
	\arrow[from=1-1, to=2-3]
	\arrow[from=1-1, to=3-2]
	\arrow[from=2-2, to=2-3]
	\arrow[from=2-2, to=3-2]
	\arrow[from=2-3, to=3-3]
	\arrow[from=3-2, to=3-3]
\end{tikzcd}$$
for any object $S \in \mathscr{S}$ and maps $S \ar X(R_0), S\ar X(R_1).$ Note here that $R_i \ar R_{01}$ are square zero (their kernel is order 2 nilpotent) extensions of $R,$ i.e. surjections $\pi_*R_i \ar \pi_*R_{01}. $ The stamenent of 2 is equivalent to saying the diagram over the image of $X$ is a \textit{pullback square}.
\end{enumerate}

\subsubsection{Remark/Explication}
Let $\text{CAlg}^{cn}_\C$ be the $\infty$-category of connective ($\pi_iR = 0$ for $i < 0$) $\mathbb{E}_\infty$-algebras over $\C.$ Let $X: \text{CAlg}^{cn}_\C \ar \mathscr{S} $ be a functor. Given a point $x \in X(\C),$ we define the formal completion of $X$ at the point $x$ to be the functor $X^\wedge: \text{CAlg}^{sm}_\C \ar \mathscr{S}$ given by 
$$X^\wedge(R) = X(R) \times_{X(\C)}\{x\}.$$
Note that the space $X^\wedge(\C)$ is automatically contractible. However, condition 2 from Definition 2.1.6 is more obtuse. We will try to illustrate with a general example. Suppose there exists some $\infty$-category $\msc$ of algebro-geometric objects such that we can do two things. Firstly, to any $A \in \text{CAlg}^{cn}_\C$, we can assign an object $\Spec A \in \msc$ which is contravariantly functorial in $A.$ Secondly, suppose there exists a special object $\mathscr{X} \in \msc$ such that $\mathscr{X}$ represents the functor $X.$ In other words, we have 
$$X(A) \simeq \Hom_\msc(\Spec A, \mathscr{X})$$
for any small $\C$-algebra $A.$ In order to verify condition 2 in this context, we can show that when $\phi: R_0 \ar R_{01}$ and $\phi': R_1 \ar R_{01}$ induce surjections $\pi_0 \ar \pi_0R_{01} \longleftarrow \pi_0R_1,$ the diagram 
$$\begin{tikzcd}
	{\Spec R_{01}} & {\Spec R_1} \\
	{\Spec R_0} & {\Spec(R_1\times_{R_{01}}R_0)}
	\arrow[from=1-1, to=1-2]
	\arrow[from=1-1, to=2-1]
	\arrow[from=1-2, to=2-2]
	\arrow[from=2-1, to=2-2]
\end{tikzcd}$$
is a pushout square ($\Spec(R_1\times_{R_{01}}R_0)$ is the colimit of the diagram 
$$\begin{tikzcd}
	{\Spec R_{01}} & {\Spec R_1} \\
	{\Spec R_0} & {}
	\arrow[from=1-1, to=1-2]
	\arrow[from=1-1, to=2-1]
\end{tikzcd}\Bigg).$$

\subsubsection{A word of warning (not to be interpreted as foreboding)}
In general, if $X$ is a formal moduli problem over $\C,$ one can always restrict $X$ to the subcategory of $\text{CAlg}^{sm}_\C$ consisting of the ordinary local artinian algebras (i.e.,$N(\text{Ring}^{art}_\C)$ to obtain a classical formal moduli problem $X_0$ with values in $\mathscr{S}.$ However the converse is not necessarily true. If we are given a classical formal moduli problem $X_0$, there need not exist a formal moduli problem $X$ with $X|_{N(\text{Ring}^{art}_\C)} = X_0.$ An example where this is in fact true is the one outlined in Remark/Explication 2.1.7.

\section{Differential graded Lie algebras and their (co)homology}

\subsection{} In this section we will introduce some terminology and constructions surrounding the concept of a differential graded (or, sometimes, dg) Lie algebra, its homology and cohomology, and begin the see some of the connections with formal moduli problems.

\subsubsection{Definition} Let $k$ be a field. A \textit{differential graded Lie algebra} $\g_*$ over $k$ is a $\Z$-graded vector space 
$$\g_* = \bigoplus \g_i$$
equipped with a differential map
$$d: \g_i \ar \g_{i-1}, \text{  } d^2 = 0$$ and a Lie bracket 
$$[-,-]: \g_p \otimes_k \g_q \ar \g_{p+q}$$
 satisfying 
 $$[x_p, x_q] + (-1)^{pq}[x_q,x_p] = 0$$
 where $x_p,x_q \in \g_p,g_q$ respectively. Moreover, if we let $x_\ell \in \g_\ell,$ the bracket satisfies the graded Jacobi identity, that is 
 $$(-1)^{p\ell}[x_p,[x_q,x_\ell]] + (-1)^{pq}[x_q, [x_\ell, x_p]] + (1)^q\ell[x_\ell[x_p,x_q]] = 0.$$
 We also maintain that $d$ is a derivation with respect to the bracket. We view  $(\g_*,d)$ as a chain complex, which gives us the impulse to make note of some categorical notions. 

 \subsubsection{A forgetful subsubsection} We denote by $\text{Vect}^{dg}_k$ the category of differential graded vector spaces over a field $k.$ The objects in this category are chain complexes 
 $$\cdots \ar V_1 \ar V_0 \ar V_{-1} \ar \cdots$$
 Note that $\text{Vect}^{dg}_k$ is a symmetric monoidal category, with the tensor product structure given by 
 $$(V \otimes W)_n = \bigoplus_{p+q = n} V_p \otimes_k W_q $$
 and the symmetric isomorphism
 $$V \otimes W \simeq W \otimes V$$
 being 
 $$\bigoplus_{p+q=n}V_p \otimes_k W_q \simeq \bigoplus_{p+q=n}W_q \otimes_k V_p$$
 multiplied by the factor $(-1)^{pq}.$ Let $V$ be a graded vector space over $k,$ and let $V^{\vee}$ be its graded dual such that 
 $$(V^\vee)_p = \Hom_k(V_{-p}, k).$$
 For all $n \in \Z,$ let $V[n]$ denote the graded shift of $V$ by $n;$ hence $V[n]_p = V_{p-n}.$
 Let $\text{Alg}^{dg}_k$ denote the category of differential graded associative $k$-algebras (more precisely, associative algebra objects of $\text{Vect}^{dg}_k,$ and let $\text{CAlg}^{dg}_k$ denote the category of commutative associative algebra objects in $\text{Vect}^{dg}_k.$ An object $A$ in $\text{Alg}^{dg}_k$ is a chain complex $(A_*,d)$ with unit in $A_0$, and the differential $d$ satisfying 
 $$d(x_p,x_q) = dx_px_1 + (-1)^pxdy$$
 where $x_p \in A_p, x_q \in A_q.$ We say that $A$ is commutative if  
 $$x_px_q = (-1)^{pq}x_qx_p.$$
 Finally, we define $\text{Lie}^{dg}_k$ to be the category of differential graded Lie algebras over $k.$ The morphisms in this category are morphisms of chain complexes which respect the Lie bracket. That is, a morphism $\varphi: (\g_*, d) \ar (\g'_*, d')$ satisfies 
 $$\varphi([x_p,x_q]) = [\varphi(x_p), \varphi(x_q)].$$

 \subsubsection{Remark} Let $A= (A_*, d)$ be a differential graded  algebra over $k.$ Then $A_*$ has the structure of a dg-Lie algebra, by 
 $$[-,-]: A_p \otimes_k A_q \ar A_{p+q}$$
 given by 
 $$[x_p,x_q] = x_px_q - (-1)^{pq}x_qx_p.$$
 This determines a forgetful functor $\text{Alg}^{dg}_k \ar \text{Lie}^{dg}_k$ with left adjoint 
 $$U: \text{Lie}^{dg}_k \ar \text{Alg}^{dg}_k$$
 which is given by assigning $g_*$ to its universal enveloping algebra, defined 
 $$\g_* \longmapsto U(\g_*) := \bigoplus_{n\ge 0}\g_*^{\otimes_n}/\Big((x\otimes y) - (-1)^{pq}(y\otimes x)-[x,y]\Big)$$
 where $x \in \g_p, y \in \g_q.$ $U(\g_*)$ in fact admits a filtration 
 $$U(\g_*)^{\le 0} \subset U(\g_*)^{\le 1} \subset \ldots $$
 where each $U(\g_*)^{\le n}$ is the image of $\bigoplus_{0\le i\le n}\g_*^{\otimes_i}$ in $U(\g_*).$

 \subsubsection{Definition} Let $\phi: \g_* \ar \g'_*$ be a morphism of dg-Lie algebras over $k.$ We say that $\phi$ is a \textit{quasi-isomorphism} if the underlying map of chain complexes induces an isomorphism on homology. 

 \subsubsection*{} We now come to an important construction in the general theory; that of a model category. See the appendix for further exposition!

 \subsubsection{Remark} The category $\text{Vect}^{dg}_k$ has the structure of a model category, wherein we say that a map of chain complexes $f: V_* \ar W_*$ is 
 \begin{enumerate}
     \item a fibration if each induced map $V_n \ar W_n$ is surjective
     \item a cofibration if each induced map $V_n \ar W_n$ is injective
     \item a weak equivalence if it is a quasi-isomorphism.
 \end{enumerate}

 \subsubsection{Proposition}
\textit{ Let $k$ be a field of characteristic zero. The category} $\text{Lie}^{dg}_k$ \textit{has the structure of a left proper combinatorial model category}.

\subsubsection{Lemma in aid of proposition 3.1.6}
\textit{Let $f: \g_* \ar \g'_* $ be a morphism of differential graded Lie algebras over $k.$ The following are equivalent:}
\begin{enumerate}
    \item $f$ is a quasi-isomorphism
    \item The induced map $U(\g_*) \ar U(\g'_*)$ is a quasi isomorphism of differential graded algebras.
    \end{enumerate}
    \textit{Proof of lemma.} For every $n \in \Z,$ let $\psi: \g_*^{\otimes_n} \ar U(\g_*)$ denote the multiplication map. For any permutation $\sigma \in \{1,2,\ldots,n\},$ let $\phi_\sigma $ be the induced automorphism of $\g_*^{\otimes_n}.$ Then the map 
    $$\frac{1}{n}\sum_\sigma \psi \circ \phi_\sigma$$
    is invariant to precomposition with $\phi_\sigma,$ and thus factors as the composition
    $$\g_*^{\otimes_n} \ar \text{Sym}^n(\g_*) \ar U(\g_*)^{\le n} \subset U(\g_*).$$
    We observe that the composition 
    $$\text{Sym}^n(\g_*) \ar U(\g_*)^{\le n} \ar \text{gr}^n(U(\g_*))$$ coincides with the isomorphism (by PBW, see [73])
    $$\theta: \text{Sym}^*(\g_*) \ar \text{gr}(U(\g_*)).$$
    It follows that the direct sum of maps 
    $$\text{Sym}^n(\g_*) \ar U(\g_*)^{\le n}$$
    is an isomorphism of chain complexes $\text{Sym}^*(\g_*) \simeq U(\g_*).$ Moreover, if $g: V_* \ar W_*$ is a quasi-isomorphism of chain complexes of $k$-vector spaces, then $g$ necessarily induces a quasi-isomorphism $\text{Sym}^*(V_*) \simeq \text{Sym}^*(W_*).$ Hence the completed proof follows by taking note of the isomorphism $\text{Sym}^*(\g_*) \simeq U(\g_*)).$ \Laughey \newline \newline
    \textit{Proof of proposition 3.1.6}  Note that the forgetful functor $\text{Lie}^{dg}_k \ar \text{Vect}^{dg}_k $ has as a left adjoint the free Lie algebra functor, denoted $\text{Free}: \text{Vect}^{dg}_k \ar \text{Lie}^{dg}_k.$ For all $n \in \Z,$ we call $E(n)_*$ the acyclic chain complex 
    $$\cdots \ar 0 \ar 0 \ar k \ar k \ar 0 \ar 0 \ar \cdots$$
    which is nontrivial only in degrees $n$ and $n-1.$ Let $\partial E(n)_*$ be the subcomplex of $E(n)_*$ only nontrivial in degree $n-1.$ Let $C_o$ be the collection of morphisms of $\text{Lie}^{dg}_k$ of the form 
    $$\text{Free}(\partial E(n)_*) \ar \text{Free}(E(n)_*),$$
    and let $W$ be the collection of quasi-isomorphisms of $\text{Lie}^{dg}_k.$ The claim here is that 
    \begin{enumerate}
        \item $W$ is \textit{perfect} (see appendix B) (this follows from Lurie [9])
        \item If $f: \g_* \ar \g'_*$ is a quasi-isomorphism over $k,$ and $x \in \g_{n-1}$ is a cycle which classifies the map 
        $$\text{Free}(\partial E(n)_*) \ar \g_k,$$
        then the induced map 
        $$\g_* \coprod_{\text{Free}(\partial E(n)_*)} \text{Free}(E(n)_*) \ar \g'_*\coprod_{\text{Free}(\partial E(n)_*)} \text{Free}(E(n)_*)$$
        is a quasi-isomorphism. 
    \end{enumerate}
    We will now sketch the proof of 2. Let $F: U(\g_*) \ar U(\g'_*)$ be the map induced by $f.$ By lemma 3.1.7, $F$ itself is a quasi-isomorphism. We can construct a differential graded algebra $B_*$ by adjoining, for the same cycle $x,$ the class $\langle y \st \deg(y) = n, dy = x\rangle$ to $U(\g_*),$ and the same for $U(\g'_*).$ We essentially want to descend to a quasi isomorphism $B_* \ar B'_*.$ We observe that $B_* $ (and respectively $B'_*)$ admits a filtration 
    $$U(\g_*) \simeq B_*^{\le 0} \subset B_*^{\le 1} \subset \ldots $$
    where each $B^{\le m}_*$ is the subspace spanned by all expressions of the form $a_0ya_1y\cdots ya_k, k \le m$ where $a_i$ are in the image of $U(\g_*)$ in $B_*.$ Since the collection of quasi-isomorphisms is perfect, it is stable under filtered colimits (appendix B), hence it suffices to show that for all $m\ge 0,$ the map 
    $$B^{\le m}_* \ar B'^{\le m}_*$$
    is a quasi-isomorphism. We will do this by induction on $m.$ The base case $m = 0$ is satisfied by assumption under the quasi-isomorphism $U(\g_*) \ar U(\g'_*).$ Let $m > 0.$ We have a diagram of short exact sequences
    $$\begin{tikzcd}
	0 & {B_*^{\le m-1}} & {B_*^{\le m}} & {B_*^{\le m}/B_*^{\le m-1}} & 0 \\
	0 & ({B'^{\le m-1}}_*) & B'^{\le m}_* & B'^{\le m}_*/B'^{\le m-1}_* & 0
	\arrow[from=1-1, to=1-2]
	\arrow[from=1-2, to=1-3]
	\arrow[from=1-2, to=2-2]
	\arrow[from=1-3, to=1-4]
	\arrow[from=1-3, to=2-3]
	\arrow[from=1-4, to=1-5]
	\arrow["\phi", from=1-4, to=2-4]
	\arrow[from=2-1, to=2-2]
	\arrow[from=2-2, to=2-3]
	\arrow[from=2-3, to=2-4]
	\arrow[from=2-4, to=2-5]
\end{tikzcd}$$
Our inductive hypothesis says that the map 
$$B^{\le m-1}_* \ar B'^{\le m-1}_* $$
is a quasi-isomorphism, so it now reduces to showing that 
$$\phi: B^{\le m}_*/B^{\le m-1}_* \ar B'^{\le m}_*/B'^{\le m-1}_* $$
is a quasi-isomorphism. We observe that the construction $a_0 \otimes \cdots \otimes a_n \longmapsto a_0ya_1y\cdots ya_m$ determines an isomorphism of chain complexes
$$U(g_*)^{\otimes_{m+1}} \ar {B_*^{\le m}/B_*^{\le m-1}}$$
(respectively for $B_*').$ This corresponds to the map 
$$U(\g_*)^{\otimes_{m+1}}\ar U(\g'_*)^{\otimes_{m+1}} $$
given by the $m+1$st tensor power of the quasi-isomorphism $F.$
\newline\newline Next, let $f: \g_* \ar \g'_*$ be a map of dg-Lie algebras with the right lifting property with respect to all morphisms in $C_0.$ We claim that $f$ is a quasi-isomorphism. Our goal is to show that $f$ induces an isomorphism $\theta_n: H_n(\g_*) \ar H_n(\g'_*)$ of the homology groups of the underlying chain complexes. First we show surjectivity. Let $\eta\in H_n(\g'_*)$ be a class represented by a cycle $x \in \g'_n.$ Then $x$ determines a map 
$$u: \text{Free}(E(n)_*) \ar \g'_*$$ which vanishes on $\text{Free}(\partial E(n)_*).$ Let $v: \text{Free}(E(n)_*) \ar \g_*$ be a map of dg-Lie algebras vanishing on $\text{Free}(\partial E(n-1)_*).$ Then $u = f\circ v$ and $v$ is determined by a cycle $\Bar{x} \in \g_n,$ which represents a homology class lifting $\eta.$ Now we wish to show injectivity.  Suppose $\eta \in H_n(\g_*)$ is a class whose image in $H_n(\g'_*)$ vanishes. Then $\eta$ is represented by a cycle $x \in \g_n$ such that $f(x) = dy$ for some $y \in \g_{n+1}.$ So $y$ determines a map $u: \text{Free}(E(n+1)_*) \ar \g'_*$ such that $u$ restricted to $\text{Free}(E(n+1)_*)$ lifts to $\g_*.$ Thus $u = f\circ v$ for some 
$$v: \text{Free}(E(n+1 )_*) \ar \g_*$$ 
whose restriction to $\text{Free}(\partial E(n+1)_*)$ classifies the cycle $x.$ Thus $x$ is a boundary, so $\eta = 0.$
It follows (see T.A.2.6.13) that $\text{Lie}^{dg}_k$ has the structure of a left proper combinatorial model category with $W$ being the class of weak equivalences, and $C_0$ the generating cofibrations. To wrap up the proof, we just need to show that a morphism $\varphi: \g_* \ar \g'_*$ in $\text{Lie}^{dg}_k$ is a fibration if and only if it is degreewise surjective. We can do this by first recognizing that if $\varphi$ is a fibration, the map of dg-Lie algebras associated to $i_n: 0 \ar \text{Free}(E(n)_*)$ factors as 
$$0 \ar 0 \coprod_{\text{Free}(\partial E(n-1)_*)}\text{Free}(E(n-1)_*)\simeq \text{Free}(\partial E(n)_*) \ar \text{Free}(E(n)_*)$$
and is thus a cofibration. Note that $E(n)$ is acyclic and so each tensor power $E(n)^{\otimes_m}_*, m>0$ is itself acyclic, hence the map 
$$k \simeq U(0) \ar U(\text{Free}(E(n)_*))\simeq \bigoplus_{m\ge 0}E(n)^{\otimes_m}_*.$$
Hence $i_n$ is a trivial cofibration such that $\varphi$ has right lifting with respect to $i_n.$ Thus $\g_n \ar \g'_n$ is surjective. Conversely, suppose that $\varphi$ is degreewise surjective. Let $S$ be the collection of all trivial cofibrations in $\text{Lie}^{dg}_k$ with left lifting (say it ten times fast) with respect to $\varphi.$ Let $f: \h_* \ar \h''_*$ be a trivial cofibration in $\text{Lie}^{dg}_k.$ We'll show that $f \in S.$ We can factor $f$ as the composition
$$\begin{tikzcd}
	{\h_*} & {\h'_*} & {\h''_*}
	\arrow["{f'}", from=1-1, to=1-2]
	\arrow["{f''}", from=1-2, to=1-3]
\end{tikzcd}$$
where $f'\in S$ and $f''$ has right lifting for each $i_n$ which $f$ contains. In otherwords, $f''$ is degreewise surjective. $f$ and $f'$ are quasi-isomorphisms, i.e. $f'$ and $f''\circ f'$ are, which implies that $f''$ is a quasi-isomorphism as well. It follows that $f''$ is a trivial fibration in the category of chain complexes and therefore is the same in $\text{Lie}^{dg}_k.$ Since $f$ is a cofibration, the diagram 
$$\begin{tikzcd}
	{\h_*} & {\h'_*} \\
	{\g''_*} & {\g''_*}
	\arrow["{f'}", from=1-1, to=1-2]
	\arrow[from=1-1, to=2-1]
	\arrow["{f''}", from=1-2, to=2-2]
	\arrow[dashed, from=2-1, to=1-2]
	\arrow[equal, from=2-1, to=2-2]
\end{tikzcd}$$
admits a completion $\g''_*\ar \h'_*.$ Thus $f$ is a retract of $f'$ and therefore $f \in S,$ which completes the proof. \Laughey 

\subsubsection{Remark} The forgetful functor $\text{Alg}^{dg}_k \ar \text{Lie}^{dg}_k$ preserves fibrations and weak equivalences, and is as such a right Quillen functor. Moreover, the universal enveloping algebra $U: \text{Lie}^{dg}_k \ar \text{Alg}^{dg}_k $ is a left Quillen functor.

\subsubsection{Proposition}
\textit{Let $\mathscr{J}$ be a small category such that $N(\mathscr{J})$ is sifted. The forgetful functor }
$$G: \text{Lie}^{dg}_k \ar \text{Vect}^{dg}_k$$
\textit{preserves $\mathscr{J}$-indexed homotopy colimits}

\subsection{Interlude into some emergent symbioses }
We wish to begin leading ourselves into the connection between formal moduli problems and differential graded Lie algebras. For starters, we have a powerful theorem connecting the two, the proof of which can be found in Pridham's [15].

\subsubsection{Theorem (powerful)}
\textit{Let} Moduli \textit{denote the full subcategory of} $\text{Fun}(\text{CAlg}^{sm}_\C, \mathscr{S})$ \textit{spanned by all of the formal moduli problems. Then there is a functor}
$$\theta: N(\text{Lie}^{dg}_\C) \ar \text{Moduli}$$
\textit{with the universal property that for every $\infty$-category $\msc,$ composition with $\theta$ induces a fully faithful embedding} $\text{Fun}(\text{Moduli}, \msc) \ar \text{Fun}(N(\text{Lie}^{dg}_\C), \msc)$ \textit{whose essential image is the collection of all functors} $F: N(\text{Lie}^{dg}_\C) \ar \msc$ \textit{carrying quasi-isomorphisms of differential graded Lie algebras to equivalences in $\msc.$}

\subsubsection{Remark}
To demonstrate the validity of this theorem's namesake, one is invited to contemplate the following. Let $W$ be the collection of quasi-isomorphisms of $\text{Lie}^{dg}_\C.$ Let $\text{Lie}^{dg}_\C[W_{-1}]$ be the $\infty$-category obtained from the nerve $N(\text{Lie}^{dg}_\C)$ by inverting all elements of $W.$ Then the above theorem implies an equivalence of $\infty$-categories
$\text{Lie}^{dg}_\C[W_{-1}] \simeq \text{Moduli}.$ In particular, every differential graded Lie algebra $\g_*$ determines a formal moduli problem. This significant result is the cornerstone of the connection between Lie algebras and commutative algebras, and finds itself one of the main focuses of this paper. This connection is controlled by the Chevalley-Eilenberg functor (or the Lie algebra cohomology of a dg-Lie algebra $\g_*),$ which assigns $\g_*$ to the cochain complex of vector spaces $C^*(\g_*).$ Particularly, this construction determines a functor 
$$(\text{Lie}^{dg}_\C)^{\text{op}} \ar \text{CAlg}^{dg}_\C,$$
carrying quasi-isomorphisms to quasi-isomorphisms, and in so doing induces a functor between $\infty$-categories 
$$(\text{Lie}^{dg}_\C[W_{-1}])^{\text{op}} \ar \text{CAlg}^{dg}_\C[W'^{-1}],$$
where $W'$ is the $\text{CAlg}^{dg}_\C$ version of $W.$ 

\subsubsection*{}
Now that we understand a lot of the machinery of differential graded Lie algebras, we can expand our categorical understanding of $\text{Lie}^{dg}_k.$ Recall that, given a model category $\msc,$ we can get its homotopy category $h\msc$ by effectively inverting its weak equivalences. The homotopy category gives us a nice way of condensing information about composable morphisms into just $\pi_0$, which throws out a lot of the extra higher order information that we have in the ambient model category. But this also might kill too much; if we care about higher order information, we might get equivalences in the homotopy category that fail in higher $\pi_n$ for $n > 0.$ So we have a bit of a goldilocks-style dilemma, and for this we turn to the notion of an underlying $\infty$-category. Let's examine our category of interest, $\text{Lie}^{dg}_k.$ Let $k$ be a field of characteristic zero. We define the underlying $\infty$-category $\text{Lie}_k$ of $\text{Lie}^{dg}_k$ to be an $\infty$-category equipped with a functor 
$$u: N(\text{Lie}^{dg}_k) \ar \text{Lie}_k$$
satisfying the (now familiar) universal property that for any $\infty$-category $\msc,$ composition with $u$ induces an equivalence from $\text{Fun}(\text{Lie}_k, \msc)$ to the full subcategory of $\text{Fun}(N(\text{Lie}^{dg}_k), \msc)$ spanned by functors $F: \text{Lie}^{dg}_k \ar \msc$ which carry quasi-isomorphisms to equivalences in $\msc.$ (By equivalences, we mean isomorphisms up to higher homotopy). We call $\text{Lie}_k$ the \textit{$\infty$-category of differential graded Lie algebras over $k.$}

\subsection{Homology and cohomology of differential graded Lie algebras}
Let $\g$ be a Lie algebra over a field $k,$ and let $U(\g)$ be its universal enveloping algebra. We can view $k$ as a left or right $U(\g)$-module where each $x \in \g$ acts trivially on $k.$ We define the homology and cohomology groups of $\g$ to be 
$$H_n(\g) = \text{Tor}^{U(\g)}_n(k,k), \text{  } H^n(\g) = \text{Ext}^n_{U(\g)}(k,k).$$
We will shortly exposit a more precise definition of the (co)homology groups, centered around the construction of the (co)homology of the \textit{Chevalley-Eilenberg complexes}. But first we review another important construction back in the setting of differential graded Lie algebras:

\subsubsection{Definition: the cone on \texorpdfstring{$\g_*$}{\g_*}}
 Let $\g_*$ be a dg-Lie algebra. We define the \textit{cone on} $\g_*,$ denoted $Cn(\g)_*$, to be a differential graded Lie algebra given by:
 \begin{enumerate}
     \item for all $n \in \Z,$ we define the vector space $Cn(\g)_*$ by $\g_n \oplus \g_{n-1}.$ Elements of $Cn(g)_n$ are of the form $x + \eps y,$ where $x \in \g_n, y \in \g_{n-1}.$
     \item The differential satisfies $d(x+\eps y) = dx+y -\eps dy$
     \item The bracket is given by $[x+\eps y, x'+\eps y'] = [x,x'] + \eps([y,x']+ (-1)^p[x,y'])$ for $x \in \g_p.$
 \end{enumerate}

\subsubsection{The homological Chevalley-Eilenberg complex}
Let $\g_*$ be a differential graded Lie algebra over a field $k.$ The zero map $\g_* \ar 0$ sneakily induces a map of differential graded algebras $U(\g)_* \ar U(0) \simeq k$. Hence there is a map of dg-Lie algebras $\g_* \ar Cn(\g_*).$ We define the \textit{homological Chevalley Eilenberg complex of }$\g_*$ to be the chain complex given by the tensor 
$$C_*(\g_*) := U(Cn(\g)_*)\otimes_{U(\g_*)} k$$

\subsubsection{Remark} We can regard the shifted chain complex $\g_*[1]$ as an abelian graded Lie algebra, and so we have a map $\g_*[1] \ar Cn(\g)_*$ (note that this is not a map of differential graded Lie algebras, so there is no differential here) inducing a map 
$$\text{Sym}^*(\g_*[1]) \ar U(Cn(\g)_*)$$
of graded vector spaces, under the identification $\text{Sym}^*(\g_*[1]) \simeq U(\g_*[1]).$ By Poincare-Birkhoff-Witt, this gives an isomorphism 
$$U(Cn(\g)_*) \simeq \text{Sym}^*(\g_*[1]) \otimes_k U(\g_*)$$
of graded right $U(\g_*)$-modules, and hence an isomorphism of graded vector spaces $\phi: \text{Sym}^*(\g_*[1]) \ar C_*(\g_*).$
Identifying $C_*(\g_*)$ with $\text{Sym}^*(\g_*[1])$ under the map $\phi,$ the differential on $C_*(\g_*)$ is given by 
$$D(x_1,\ldots,x_n) = \sum_{1\le i\le n}(-1)^{p_1+\cdots+p_{i-1}}x_1\cdots dx_i\cdots x_n$$
$$+ \sum_{1\le i<j\le n}(-1)^{p_i(p_{i+1}+\cdots p_{j-1})}x_1\cdots \hat{x_i}\cdots \hat{x_j}\cdots x_n[x_i,x_j].$$

\subsubsection{Remark}
The filtration of $\text{Sym}^*(\g_*)$ by $\bigoplus_{i\le n}\text{Sym}^i(\g_*)$ defines a filtration 
$$k\simeq C^{\le 0}_*(\g_*) \hookrightarrow C^{\le 1}_*(\g_*) \hookrightarrow C^{\le 2}_*(\g_*) \hookrightarrow \cdots.$$
Moreover, using the formula for $D(x_1,\ldots,x_n),$ we obtain the canonical isomorphisms 
$$C^{\le n}_*(\g_*)/C^{\le n-1}_*(\g_*) \simeq \text{Sym}^n(\g_*)$$
of differential graded $k$-vector spaces. 

\subsubsection{Proposition} 
\textit{Let $k$ be a field of characteristic zero, and let $f: \g_* \ar \g'_*$ be a quasi-isomorphism of dg-Lie algebras. Then the induced map on the homological CE-complexes $C_*(\g_*) \ar C_*(\g'_*)$ is a quasi-isomorphism of chain complexes.}
\newline\newline 
\textit{Proof.} Since the collection of quasi-isomorphisms is closed under filtered colimits, it suffices to show that the map 
$$\theta_n: C^{\le n}_*(\g_*) \ar C^{\le n}_*(\g'_*)$$
is a quasi-isomorphism for each $n.$ We proceed by induction on $n.$ If $n=0$ there is an immediate isomorphism, so we assume $n> 0.$ We have a commutative diagram of short exact sequences
$$\begin{tikzcd}
	0 & {C^{\le n-1}_*(\g_*)} & {C^{\le n}_*(\g_*)} & {\text{Sym}^n(\g_*[1])} & 0 \\
	0 & {C^{\le n-1}_*(\g'_*)} & {C^{\le n}_*(\g'_*)} & {\text{Sym}^n(\g'_*[1])} & 0
	\arrow[from=1-1, to=1-2]
	\arrow[from=1-2, to=1-3]
	\arrow[from=1-2, to=2-2]
	\arrow[from=1-3, to=1-4]
	\arrow[from=1-3, to=2-3]
	\arrow[from=1-4, to=1-5]
	\arrow["\phi", from=1-4, to=2-4]
	\arrow[from=2-1, to=2-2]
	\arrow[from=2-2, to=2-3]
	\arrow[from=2-3, to=2-4]
	\arrow[from=2-4, to=2-5]
\end{tikzcd}$$
By our inductive hypothesis, it suffices to show that the map $\phi$ between the symmetric algebras is a quasi-isomorphism. And since char $k = 0,$ $\phi$ is a retract of the map $\g^{\otimes_n}_*[7] \ar \g'^{\otimes_n}_*[7]$ which is a quasi-isomorphism by the assumption that $f$ is. \Laughey
\newline\newline Looking forward, if $\g_*$ is a differential graded Lie algebra, we call the homology groups of $C_*(\g_*)$ the \textit{Lie algebra homology groups of} $\g_*.$

\subsubsection{The cohomological Chevalley-Eilenberg complex}
Let $\g_*$ be a differential graded Lie algebra over a field $k.$ We denote by $C^*(\g_*)$ the linear dual of $C_*(\g_*),$ which we call the \textit{cohomological Chevalley-Eilenberg complex.} Much of the work constructing $C^*(\g_*)$ has already been done, but it is significant to examine the natural multiplication structure of $C^*(\g_*),$ carrying $\lambda \in C^p(\g_*) $ and $\mu \in C^q(\g_*)$ to $C^{p+q}(\g_*).$ We identify elements of $C^p(\g_*)$ with the dual space of the graded vector space $\text{Sym}^p(\g_*[1]).$ Let 
$$S = \{i_1 < \ldots < i_m\}, S' = \{j_1 < \ldots < j_{n-m}\},$$ so $S \cup S' = \{1,\ldots,n\}.$ Write $p = r_{i_1} + \ldots + r_{i_m}.$ We define, for $x_i \in \g_{r_i},$ 
$$(\lambda\mu)(x_1,\ldots,x_n) = \sum \prod_{i \in S,j\in S', i<j}(-1)^{r_ir_j}\lambda(x_{i_1}\ldots x_{i_m})\mu(x_{j_1}\ldots x_{j_{n-m}}).$$
With this multiplication, $C^*(\g_*)$ has the structure of a commutative differential graded algebra.

\section{Weaving together}
\subsection{}Recall our proposition of great importance. As it turns out, the PGI has a converse, which we will denote by coPGI, which stipulates that a formal moduli problem $X$ is determined by $\g_*$ up to equivalence. More precisely, we would like to prove the following co-proposition of great importance: 
\subsubsection{Theorem (coPGI)}
\textit{Let $k$ be a field of characteristic zero. Let} $\text{Lie}_k$ \textit{denote the $\infty$-category underlying} $\text{Lie}^{dg}_k$ \textit{obtained by inverting quasi-isomorphisms. Then there is an equivalence of $\infty$-categories} 
$$\Psi: \text{Lie}_k \ar \text{Moduli}_k.$$

\subsection{Koszul Duality}
Let $k$ be a field of characteristic zero. It follows from 3.3.5 that the functor $C^*: \g_* \longmapsto C^*(\g_*)$ carries quasi-isomorphisms to quasi-isomorphism. We then obtain a functor between $\infty$-categories $\text{Lie}_k \ar \text{CAlg}^{op}_k,$ which we still denote by $C^*.$ Note that this functor carries the intial object $0 \in \text{Lie}_k$ to the terminal object $k \in \text{CAlg}^{op}_k.$ We obtain another functor $\text{Lie}_k \ar (\text{CAlg}^{aug}_k)^{op}$ whose target is the $\infty$-category of augmented $\mathbb{E}_\infty$-algebras over $k.$ We continue to abuse notation by calling this functor $C^*$ as well. This functor preserves small colimits, and we note that $\text{Lie}_k$ is presentable. We define the functor 
$$\mf{D}: (\text{CAlg}^{aug}_k)^{op} \ar \text{Lie}_k$$
to be the right adjoint of the functor $C^*: \text{Lie}_k \ar (\text{CAlg}^{aug}_k)^{op}, $ and call it the \textit{Koszul duality functor}. Our goal in introducing this is to prove that $\mf{D}$ is a deformation theory, which will help us prove the coPGI. First we need to verify that $\mf{D}$ is a weak deformation theory. Recall these axioms. The first two are automatic, since $\text{Lie}_k$ is presentable and $\mf{D}$ admits a left adjoint by construction. For axiom 3, we will prove the following 
\subsubsection{Proposition}
Let $k$ be a field of characteristic zero, and $\g_*$ a differential graded Lie algebra over $k.$ Let $\msc$ be the full subcategory of $\text{Lie}_k$ spanned by cofibrant (with respect to the model on $\text{Lie}^{dg}_k$) objects satisfying the following
\begin{enumerate}
    \item There exists a graded vector space $V_* \subset \g_*$ such that for each integer $n,$ $\dim V_n < \infty.$
    \item For all $n \ge 0,$ $V_n$ is trivial
    \item $V_*$ freely generates $\g_*$ as a graded Lie algebra.
\end{enumerate}
Then $\msc$ satisfies axiom 3. The proof of this relies on the following lemma, whose proof can be found in Lurie [7] \S2.
\subsubsection{Lemma in aid of Proposition}
\textit{Let $\g_*$ be a differential graded Lie algebra over $k,$ and assume that for each $n,$ $\dim \g_n < \infty,$ and that $\g_n$ is trivial for each $n\ge 0.$ Then the unit map $u: \g_* \ar \mf{D}(C^*(\g))$ is an equivalence in $\text{Lie}_k.$}
\newline\newline We are now approaching the proof of the coPGI. We just need to know how to construct the functor $\Psi: \text{Lie}_k \ar \text{Moduli}^\Gamma.$ Let $\g_* \in \text{Lie}^{dg}_k,$ and $R \in \text{CAlg}^{sm}_k.$ We can identify $R$ with an augemented commutative dg-algebra over $k.$ Call its augmentation ideal $\mf{m}_R.$ Then the tensor product $\mf{m}_R \otimes_k \g_*$ is a differential graded Lie algebra over $k.$ To properly construct $\Psi,$ we want $\Psi(\g_*)(R)$ to be a space of \textit{Maurer-Cartan elements}, i.e the space of solutions to the Maurer-Cartan equation $dx = [x,x].$ We call such a space $MC(\g_*).$ Fortunately for us, there is a well defined bifunctor 
$$MC: \text{CAlg}^{aug}_k \times \text{Lie}_k \ar \mathscr{S}$$
given by $(R, \g_*) \longmapsto MC(\mf{m}_R\otimes_k \g_*),$ which we can also describe in terms of the Koszul duality functor, namely 
$$MC(R,\g_*) = \text{Map}_{\text{Lie}_k}(\mf{D}(R), \g_*).$$
This is how we'll define our functor $\Psi.$

\subsubsection{Theorem [7]}
\textit{Let $k$ be a field of characteristic zero. Let} $(\text{CAlg}^{aug}_k, \{E\})$ \textit{be the deformation context we work with. Then the Koszul duality functor}
$$\mf{D}: (\text{CAlg}^{aug}_k)^{op} \ar \text{Lie}_k$$
\textit{is a deformation theory}.
\newline\newline \textit{Sketch of proof.}
Let $E \in \text{Stab}(\text{CAlg}^{aug}_k)$ be the spectrum object corresponding to $k,$ such that $\Omega^{\infty-n}E \simeq $ the square zero extension $k \oplus k[n].$ The previous proposition shows that $\mf{D}(E)$ is given by the infinite loop object $\{\text{Free}(k[-n-1])\}_{n\ge 0}$ in $\text{Lie}^{op}_k$ (see Lurie [7]). Here Free$ :\text{Mod}_k \ar \text{Lie}_k$ denotes the left adjoint of the forgetful functor $\theta: \text{Lie}_k \ar \text{Mod}_k.$ It follows that the functor $e: \text{Lie}_k \ar \textbf{Spc} $ from the prior chapter is given by the composition $(F \circ \theta)[1],$ where $F: \text{Mod}_k = \text{Mod}_k(\textbf{Spc}) \ar \textbf{Spc}$ is the forgetful functor. The claim of this proof follows from the technical $\infty$-categorical considerations of observing that $F$ and $\theta$ preserve colimits and sifted colimits, respectively. \Laughey
\newline\newline We are now ready to prove the main result: 

\subsubsection{Theorem (Lurie) a.k.a coPGI}
\textit{Let $k$ be a field of characteristic zero. Let} $\text{Lie}_k$ \textit{denote the $\infty$-category underlying} $\text{Lie}^{dg}_k$ \textit{obtained by inverting quasi-isomorphisms. Then there is an equivalence of $\infty$-categories} 
$$\Psi: \text{Lie}_k \ar \text{Moduli}_k.$$
\newline\newline\textit{Proof.} Let $k$ be a field of characteristic zero. Let $\Psi: \text{Lie}_k \ar \text{Fun}(\text{CAlg}^{sm}_k, \mathscr{S})$ denote the functor given by objects of the form
$$\Psi(\g_*)(R) = \text{Map}_{\text{Lie}_k}(\mf{D}(R), \g_*).$$ Combining theorems 1.4.6 and 4.2.3, we observe that $\Psi$ is a fully faithful embedding whose essential image (smallest subcategory respecting isomorphism which contains the image) is the fullsubcategory $\text{Moduli}_k \subset \text{Fun}(\text{CAlg}^{sm}_k, \mathscr{S})$ spanned by formal moduli problems. \Laughey

\newpage

\section{Appendix A: Dan Quillen and model categories}

\subsection*{}Model categories were introduced by Dan Quillen in his book \textit{Homotopical Algebra} from 1967, for the purpose of providing a framework for homotopy theory A \textit{model category} is roughly speaking a category $\mathscr{C}$ which has three classes of morphisms, called \textit{weak equivalences, fibrations,} and \textit{cofibrations.} Weak equivalences play (the more generalized) role of homotopy equivalences, while fibrations and cofibrations are more like inclusions and surjections, respectively, satisfying some lifting properties. This is of course very abstract and un-concrete, but we hope to resolve this with some examples, and, more importantly, the connection to the main subject matter of this paper. 

\subsection{Model categories}

 Suppose that $\mathscr{C}$ is a category. A morphism $f \in \mathscr{C}$ is called a \textit{retract} of a map $g \in \msc$ if there exists a commutative diagram of the form
$$\begin{tikzcd}
	A & C & A \\
	B & D & B
	\arrow[from=1-1, to=1-2]
	\arrow["f"', from=1-1, to=2-1]
	\arrow[from=1-2, to=1-3]
	\arrow["g"', from=1-2, to=2-2]
	\arrow["f", from=1-3, to=2-3]
	\arrow[from=2-1, to=2-2]
	\arrow[from=2-2, to=2-3]
\end{tikzcd}$$
such that compositions $A \ar C \ar A = \text{id}_A$ and $B \ar D \ar B = \text{id}_B.$ A (functorial) factorization is a pair of functors $(\alpha, \beta)$ from $\text{Map}(\msc) \ar \text{Map}(\msc)$ such that $f = \beta(f) \circ \alpha(f)$ with agreements 
$$\begin{tikzcd}
	{\text{domain of }f} & {\text{domain of } \alpha(f)} & {\text{domain of }\beta(f)} \\
	{\text{codomain of }f} & {\text{codomain of }\beta(f)} & {\text{codomain of }\alpha(f)}
	\arrow[equal,  from=1-2, to=1-1]
	\arrow[equal, from=2-1, to=2-2]
	\arrow[equal, from=2-3, to=1-3]
\end{tikzcd}$$
Suppose that $\psi: A \ar B$ and $\varphi: X \ar Y$ are maps in $\msc.$ We say that $\psi$ has the \textit{left lifting property} with respect to $\varphi$ and that $\varphi$ has the \textit{right lifting property} with respect to $\psi$ if for every commutative square 
$$\begin{tikzcd}
	A & X \\
	B & Y
	\arrow["f", from=1-1, to=1-2]
	\arrow["\psi"', from=1-1, to=2-1]
	\arrow["\varphi", from=1-2, to=2-2]
	\arrow["h", dashed, from=2-1, to=1-2]
	\arrow["g"', from=2-1, to=2-2]
\end{tikzcd}$$
there exists a lift $h: B \ar X$ such that $h \circ \psi = f$ and $\varphi \circ h = g.$

\subsubsection{Definition}
A \textit{model structure} on a category $\msc$ is the three subcategories of $\msc$ of weak equivalences, fibrations, and cofibrations, satisfying 
\begin{enumerate}
    \item Let $f,g$ be morphisms in $\msc$ such that $g\circ f$ is definable. If any pair of the three maps $f,g,$ and $g \circ f$ are weak equivalences, then so is the third.
    \item If $f$ is a retract of $g$ and $g$ is a weak equivalence, fibration, or cofibration, then so is $f.$
    \item A map $f$ is called a trivial or acyclic (co)fibration if it is a (co)fibration and a weak equivalence. We mandate that acyclic cofibrations satisfy left lifting with respect to fibrations, and cofibrations satisfy left lifting with respect to acyclic fibrations.
\end{enumerate}
A \textit{model category} is then a category $\msc$ with a model structure and small (co)limits. By the first axiom, any model category has an initial object $\varnothing$ and a terminal object *. We call an object $X \in \msc$ \textit{fibrant} if the unique map $X \ar *$ is a fibration, and \textit{cofibrant} if the unique map $\varnothing \ar X$ is a cofibration.

\subsubsection{Quillen functors}
Suppose that $\msc$ and $\mathscr{D}$ are model categories. A \textit{left Quillen functor} $F: \msc \ar \mathscr{D}$ is a functor which is left adjoint and preserves (acyclic) cofibrations. A \textit{right Quillen functor} is similarly a functor $G: \msc \ar \mathscr{D}$ which is right adjoint and preserves (acyclic) fibrations.

\subsubsection{A pertinent type of model category}
Categories with their objects being chain complexes form important types of model categories. For instance, let $A$ be an abelian Grothendieck category. We can define a category $C(A)$ with objects being chain complexes 
$$\cdots \ar X_1 \ar X_0 \ar X_{-1} \ar \cdots$$
and morphisms being chain maps. Then $C(A)$ has a model structure by setting cofibrations as monomorphisms and weak equivalences as quasi-isomorphisms.

\subsubsection{Definition: Combinatorial model category}
A model category $\msc$ is called \textit{combinatorial} if it contains 
\begin{enumerate}
    \item a set $S$ of small objects such that every object in $\msc$ is a colimit over objects in $S.$ Equivalently, $\msc$ has a fully faithful right adjoint localization $\msc \hookrightarrow \text{Psh}(S),$ where $\text{PSh}(S)$ is the category of presheaves on $S.$
    \item a set of cofibrations and a set of acyclic cofibrations which "generate" all (acyclic) cofibrations in $\msc.$
\end{enumerate}

Note also that a model category is \textit{left proper} if for every diagram 
$$\begin{tikzcd}
	A & X \\
	B & {X \cup B}
	\arrow["f", from=1-1, to=1-2]
	\arrow["i"', hook, from=1-1, to=2-1]
	\arrow["h", from=1-2, to=2-2]
	\arrow[from=2-1, to=2-2]
\end{tikzcd}$$
where $i$ is a cofibration and $f$ is a weak equivalence, the map $h$ is also a weak equivalence.

\subsection{Quillen adjunction and equivalence}
Let $\msc$ and $\mathscr{D}$ be model categories. Suppose we have a pair of adjoint functors
$$\begin{tikzcd}
	\msc & {\mathscr{D}}
	\arrow["F", from=1-1, to=1-2]
	\arrow["G", shift left=3, from=1-2, to=1-1]
\end{tikzcd}.$$ 
\subsubsection{Proposition}
The following are equivalent: 
\begin{enumerate}
    \item $F$ preserves (trivial) cofibrations
    \item $G$ preserves (trivial) fibrations
    \item $F$ preserves cofibrations and $G$ preserves fibrations
    \item $F$ preserves trivial cofibrations and $G$ preserves trivial fibrations.
\end{enumerate}
If any of these are satisfied, we say that $(F,G)$ deetermines a \textit{Quillen adjunction} on the categories $\msc$ and $\mathscr{D}.$
\newline\newline
Suppose that 
$$\begin{tikzcd}
	\msc & {\mathscr{D}}
	\arrow["F", from=1-1, to=1-2]
	\arrow["G", shift left=3, from=1-2, to=1-1]
\end{tikzcd}$$
is a Quillen adjunction. We also say that $F$ is a left Quillen functor and $G$ is a right Quillen functor. We can get the \textit{homotopy category} $h\msc$ of $\msc$ by first passing to the full subcategory of cofibrant objects in $\msc$ and inverting all weak equivalences (similarly for $\mathscr{D}).$ Because $F$ preserves trivial cofibrations (i.e., weak equivalences between cofibrant objects), it induces a functor $h\msc\ar h\mathscr{D}$ called the \textit{left derived functor} of $F,$ denoted $LF.$ Similary, we can describe the right derived functor of $G,$ called $RG.$ Moreover, it is the case that 
$$\begin{tikzcd}
	{h\msc} & {h\mathscr{D}}
	\arrow["LF", from=1-1, to=1-2]
	\arrow["RG", shift left=3, from=1-2, to=1-1]
\end{tikzcd}$$
determines an adjunction. 
\subsubsection{Proposition}
Let
$$\begin{tikzcd}
	\msc & {\mathscr{D}}
	\arrow["F", from=1-1, to=1-2]
	\arrow["G", shift left=3, from=1-2, to=1-1]
\end{tikzcd}$$
be a Quillen adjunction. Then the following are equivalent: 
\begin{enumerate}
    \item $LF: h\msc \ar h\mathscr{D}$ is an equivalence of categories.
    \item $RG:h\mathscr{D}\ar h\msc$ is an equivalence of categories.
    \item For any cofibrant object $C \in \msc,$ and any fibrant object $D \in \mathscr{D},$ a map $C \ar G(D)$ is a weak equivalence in $\msc$ if and only if the adjoint $F(C) \ar D$ is a weak equivalence in $\mathscr{D}.$ 
\end{enumerate}
\textit{Proof.} $1. \iff 2.$ is immediate since $RG$ and $LF$ are adjoint. Both are equivalent to the statement that 
$$u: \text{id}_\msc \ar RG \circ LF, v: LF \circ RG \ar \text{id}_\mathscr{D}$$
are weak equivalences. Moreover, we have $(RG \circ LF)(C) = G(D)$ where $F(C) \ar D$ is a weak equivalence in $\mathscr{D}.$ Thus $u$ is a weak equivalence when evaluated on $C$ if and only if for every weak equivalence $F(C) \ar D,$ the adjoint $C \ar G(D)$
is a weak equivalence. The same argument follows for $v,$ and so $1 \iff 2 \iff 3.$ \Laughey \newline If any of the three hold, we sat that $(F,G)$ gives a \textit{Quillen equivalence} on $\msc, \mathscr{D}.$
\subsubsection{Homotopy limits and colimits}
Let $\msc$ be a category with weak equivalences, and let $\mathscr{D}$ be a small diagram category. We can turn the functor category $\text{Fun}(\mathscr{D}, \msc)$ (the category with objects as functors $\mathscr{D} \ar \msc$ and morphisms natural transformations) into a category with weak equivalences by declaring them to be those natural transformations which are objectwise weak equivalences. The \textit{homotopy limit} of a functor $G: \mathscr{D} \ar \msc$ is the image of $G$ under the right derived functor of the limit $\lim_\mathscr{D}: \text{Fun}(\mathscr{D}, \msc) \ar \msc$ with respect to weak equivalences on $\msc$ and  $\text{Fun}(\mathscr{D}, \msc).$ Similarly, the \textit{homotopy colimit} of a functor $H: \mathscr{D} \ar \msc$ is the image of $H$ under the left derived functor of $\text{colim}_D:\text{Fun}(\mathscr{D}, \msc) \ar \msc$ with respect to weak equivalences on $\msc$ and  $\text{Fun}(\mathscr{D}, \msc).$

\section{Appendix B: \texorpdfstring{$\infty$}{\infty}-categorical miscellany}
\subsection{"The weeds," as it were }
\subsubsection{The category \texorpdfstring{$\text{Mod}_k$}{\text{Mod}_k}}
We define the $\infty$-category $\text{Mod}_k$ to be that of $k$-module spectra, where $k$ is a field. The idea of a $k$-module spectrum is vastly denser than what the overleaf file which induced the pdf you are reading can handle; as such, we refer the reader to the nice array of references which exposit the topic in great detail. Morally speaking, one can view the objects of $\text{Mod}_k$ as being given by chain complexes of $k$-vector spaces, where, for any $M \in \text{Mod}_k,$ the homotopy groups $\pi_*M$ constitute graded $k$-vector spaces. We also say that $M$ is \textit{locally finite} if each homotopy group is finite dimensional.

\subsubsection{Augmentation}
For $A \in \text{CAlg}_k$ (denoted $\text{Alg}^{(n)}_k$ if $n \ne \infty),$ an augmentation of $A$ is a map of $\mathbb{E}_n$-algebras $\eps: A \ar k.$
\subsubsection{The \texorpdfstring{$\infty$}{\infty}-category of spaces}
We call $\mathscr{S}$ the $\infty$-\textit{category of spaces}. We can define, roughly, $\mathscr{S}$ to be the $\infty$-categorical analogue of \textbf{Set}, by replacing equalities with homotopies.
\subsubsection{Pushouts}
Let $\msc$ be an $\infty$-category. If we are given a diagram 
$$D = \begin{tikzcd}
	C & B \\
	A & P
	\arrow["{f'}", from=1-1, to=1-2]
	\arrow["f"', from=1-1, to=2-1]
	\arrow["{g'}", from=1-2, to=2-2]
	\arrow["g"', from=2-1, to=2-2]
\end{tikzcd}$$
in $\msc,$ we say that $D$ is a \textit{pushout} (or pushout \textit{square}) if for any object $X \in \msc,$ giving me a map $P \ar X$ is morally the same as giving me two maps $A \ar X, B \ar X.$ In other words, we obtain $P$ as the colimit $P \simeq A \coprod_C B$ of the diagram
$$\begin{tikzcd}
	& C \\
	A && B
	\arrow[from=1-2, to=2-1]
	\arrow[from=1-2, to=2-3]
\end{tikzcd}$$

\subsubsection{Definition: Perfection (of classes of morphisms)}
Let $\msc$ be a presentable category. A class $W$ of morphisms in $\msc$ is called \textit{perfect} (in [40] A.2.6.12) if it satisfies the following conditions: 
\begin{enumerate}
    \item Every isomorphism is an element of $W.$
    \item For any pair of composable morphisms $f$ and $g,$ if any pair of the three maps $f,g, g\circ f$ is in $W,$ then the third is as well.
    \item Let $\{f_\alpha: X_\alpha \ar Y_\alpha\}$ be a collection of morphisms indexed by a filtered poset.  Let 
    $$X = \lim_{\ar}\{X_\alpha\}, \text{  } Y = \lim_{\ar}\{Y_\alpha\}.$$
    Let $f: X \ar Y$ be the induced map. If each $f_\alpha \in W,$ then $f \in W.$ This is equivalent to the statement that $W$ is stable under filtered colimits.
    \item There exists a subset $W_0 \subset W$ such that every $f \in W$ is a filtered colimit of morphisms in $W_0.$
\end{enumerate}

\subsubsection{The nerve of a small category}
Let $\mathscr{J}$ be a small category. The \textit{nerve} $N(\mathscr{J})$ of $\mathscr{J}$ is the simplicial set whose 0-simplices are objects of $\mathscr{J},$ 1-simplices are morphisms in $\mathscr{J},$ 2-simplices are pairs of composable morphisms, so on and so forth. We say that $N(\mathscr{J})$ is \textit{sifted} if for any family of diagrams $D_1,D_2,\ldots,D_n: \mathscr{J} \ar \text{Set},$ the set-theoretic colimits of $D_i$ commute with finite products. Concretely, if $N(\mathscr{J})$ is sifted, then 
$$\colim(D_1 \times \cdots \times D_n) \simeq \colim D_1 \times \colim D_2 \times \cdots \times \colim D_n.$$

\subsubsection{The largest Kan complex}
 Let $\msc^\simeq$ denote the $\infty$-category one gets from throwing out all non invertible morphisms of $\msc.$ This is equivalently the largest Kan complex contained in $\msc.$

\subsection{Stabilizations, loops, and suspensions}
\subsubsection{The loop and suspension functors}
Let $\msc$ be a pointed $\infty$-category, i.e. $\msc$ has a 0-object and finite (co)limits. The \textit{loop functor} $\Omega: \msc \ar \msc$ takes an object $X \in \msc$ to its space of loops based at the zero map. That is, $\Omega X = 0\times_X 0$ is defined by the pullback square 
$$\begin{tikzcd}
	{\Omega X} & 0 \\
	0 & X
	\arrow[from=1-1, to=1-2]
	\arrow[from=1-1, to=2-1]
	\arrow[from=1-2, to=2-2]
	\arrow[from=2-1, to=2-2]
\end{tikzcd}$$
Suggestively and similarly, the \textit{suspension functor} $\Sigma: \msc \ar \msc$ is defined by $\Sigma X = 0\coprod_X 0,$ i.e. given by the pushout
$$\begin{tikzcd}
	X & 0 \\
	0 & {\Sigma X}
	\arrow[from=1-1, to=1-2]
	\arrow[from=1-1, to=2-1]
	\arrow[from=1-2, to=2-2]
	\arrow[from=2-1, to=2-2]
\end{tikzcd}$$

\subsubsection{Definition: Stabilization}
Let $\msc$ be a pointed $\infty$-category. The \textit{stabilization} $\text{Stab}(\msc),$ which would benefit from less conspicuous notation, is the stable $(\Omega$ and $\Sigma$ are mutual inverses) $\infty$-category of spectrum objects in $\msc.$ Less mercurially, the objects of $\text{Stab}(\msc)$ are sequences $X_0,X_1,X_2,\ldots$ with equivalences $X_n \simeq \Omega X_{n+1}.$ $\text{Stab}(\msc)$ satisfies the following universal property motivated by the existence of the a canonical functor 
$$\Sigma^\infty: \msc \ar \text{Stab}(\msc)$$
such that for any stable $\infty$-category $\mathscr{D},$ precomposition with $\Sigma^\infty$ yields an equivalence
$$\text{Fun}_{\text{ex}}(\text{Stab}(\msc), \mathscr{D}) \simeq \text{Fun}_*(\msc, \mathscr{D})$$
where the LHS consists of exact functors and the RHS of 0-preserving functors.

\newpage 
\section{References}
\begin{enumerate}
    \item  Calaque and Grivaux, \textit{Formal moduli problems and formal derived stacks}, https://arxiv.org/pdf/1802.09556
    \item nLab, "Cofibrantly generated model category", https://ncatlab.org/nlab/show/cofibrantly+generated+model+category
    \item Fulton, W \textit{Algebraic Curves}, W.A. Benjamin, Inc., New York, 1969. https://dept.math.lsa.umich.edu/~wfulton/CurveBook.pdf
    \item Greenless, J.P.C and MAy, J.P \textit{Equivariant stable homotopy theory} https://www.math.uchicago.edu/~may/PAPERS/Newthird.pdf
    \item Hartshorne, R \textit{Deformation Theory}, Springer Graduate Texts in Mathematics
    \item Hovey, M. \textit{Model categories}, https://people.math.rochester.edu/faculty/doug/otherpapers/hovey-model-cats.pdf
    \item Lurie, J. \textit{Derived Algebraic Geometry X: Formal Moduli Problems}, https://people.math.harvard.edu/~lurie/papers/DAG-X.pdf
    \item Lurie, J. \textit{Higher Topos Theory}, https://www.math.ias.edu/~lurie/papers/HTT.pdf
    \item Lurie, J. \textit{Higher Algebra}, https://www.math.ias.edu/~lurie/papers/HTT.pdf
    \item Lurie, J. \textit{Spectral Algebraic Geometry} https://www.math.uchicago.edu/~may/PAPERS/kmbooklatex.pdf
    \item Lurie, J. \textit{Derived Algebrai Geometry III: Commutative Algebra} https://math.mit.edu/~lurie/papers/DAG-III.pdf
    \item Lurie, J. \textit{Derived Algebraic Geometry IV: Deformation Theory} https://math.mit.edu/~lurie/papers/DAG-IV.pdf
    \item Riehl, E. and Verity, D. \textit{Infinite Category Theory from Scratch}, https://math.jhu.edu/~eriehl/scratch.pdf
    \item Matsumura, H \textit{Commutative Ring Theory}, Cambridge University Press, 1986.
    \item Pridham, J \textit{Unifying derived deformation theories}, Adv. Math. 224 (2010), no.3, 772-826.
    \item Lurie, J. \textit{Derived Algebraic Geometry VI: \texorpdfstring{$\mathbb{E}[k]$}{\mathbb{E}[k]}-algebras},https://people.math.harvard.edu/~lurie/papers/DAG-VI.pdf
    \item nLab, "E-infinity algebra," https://ncatlab.org/nlab/show/E-infinity+algebra
    \item Kriz, I. and May, J.P. \textit{Operads, Algebras, Modules, and Motives} https://www.math.uchicago.edu/~may/PAPERS/kmbooklatex.pdf
    \item Soergel, W \textit{Koszul Duality and applications in representation theory} https://people.mpim-bonn.mpg.de/geordie/Soergel.pdf
    
\end{enumerate}

\end{document}